\documentstyle{amsppt}

\magnification=\magstep1
\baselineskip=13pt   
\parindent=0pt
\parskip=10pt        

\vsize=7.7in
\voffset=-.4in

\overfullrule=0pt

\define\A{{\Bbb A}}
\define\C{{\Bbb C}}

\define\R{{\Bbb R}}
\define\Q{{\Bbb Q}}

\define\a{\alpha}

\define\e{\epsilon}
\redefine\l{\lambda}
\redefine\o{\omega}
\define\ph{\varphi}

\redefine\P{\Phi}
\predefine\Sec{\S}

\define\back{\backslash}
\define\lra{\longrightarrow}
\redefine\tt{\otimes}


\define\isoarrow{\ {\overset{\sim}\to{\longrightarrow}}\ }

\define\und#1{\underline{#1}}
\define\nass{\noalign{\smallskip}}
\define\triv{1\!\!1}



\predefine\oldvol{\vol}
\redefine\vol{\text{\rm vol}}

\define\CT#1{\operatornamewithlimits{CT}_{#1}}   

\font\cute=cmitt10 at 12pt
\font\smallcute=cmitt10 at 9pt
\define\kay{{\text{\cute k}}}
\define\smallkay{{\text{\smallcute k}}}

\define\tr{\text{\rm tr}}

\redefine\Re{\text{\rm Re}}

\define\Sym{\text{\rm Sym}}

\define\underrule#1{$\und{\text{#1}}$}


\define\garrett{\bf1}
\define\grosskudla{\bf2}
\define\duke{\bf3}
\define\annals{\bf4}  
\define\howeps{\bf5}
\define\ikedatriple{\bf6}
\define\ikedapoles{\bf7}
\define\IkedaDuke{\bf8}
\define\jiang{\bf9}
\define\kimshahiditriple{\bf10}
\define\kimshahidi{\bf11}
\define\annalsII{\bf12}
\define\krinvardist{\bf13}
\define\krrdps{\bf14}
\define\krannals{\bf15}
\define\psrallis{\bf16}
\define\rallisHDC{\bf17}
\define\rallisbook{\bf18}
\define\dinakar{\bf19}
\define\roberts{\bf20}
\define\sahi{\bf21}
\define\shimizu{\bf22}
\define\vtan{\bf23}
\define\zhuI{\bf24}
\define\zhuII{\bf25}

\centerline{\bf On a conjecture of Jacquet}

\centerline{\bf by}
\centerline{\bf Michael Harris
\footnote{Institut de Math\'ematiques de
Jussieu, U.M.R. 7586 du CNRS.   Membre, Institut Universitaire de France}}
\centerline{\bf and}
\centerline{\bf Stephen S. Kudla\footnote{Partially 
supported by NSF grant DMS-9970506 and by a Max-Planck Research Prize 
from the Max-Planck Society and Alexander von Humboldt Stiftung. }}

\vskip .1in
\centerline{\it For Joe Shalika, with our admiration and appreciation.}
\vskip .1in

\subheading{Introduction}

Let $\kay$ be a number field and let $\pi_i$, $i=1$, $2$, $3$ be 
cuspidal automorphic
representations of $GL_2(\A)$ such that the product of their central characters
is trivial. Jacquet then conjectured that the central value 
$L(\frac12,\pi_1\tt\pi_2\tt\pi_3)$ of the
triple product L--function is nonzero if and only if there exists a 
quaternion algebra $B$
over $\kay$ and automorphic forms $f_i^B\in \pi_i^B$ such that the integral
$$I(f_1^B,f_2^B,f_3^B) = \int_{Z(\A)B^\times(\smallkay)\back 
B^\times(\A)} f_1^B(b)f_2^B(b) f_3^B(b)\, d^\times b \ne 0,
\tag0.1$$
where $\pi_i^B$ is the representation of $B^\times(\A)$ corresponding 
to $\pi_i$.

In a previous paper \cite{\annals}, we proved this conjecture in the
special case where $\kay=\Q$ and the $\pi_i$'s correspond to a triple 
of holomorphic newforms.
Our method was based on a combination of the Garrett, 
Piatetski-Shapiro, Rallis integral representation
of the triple product L-function with the extended Siegel--Weil 
formula and the seesaw identity.
The restriction to holomorphic newforms over $\Q$ arose from (i) the 
need to invoke the Ramanujan Conjecture
to control the poles of some bad local factors and (ii) the use of a 
version of the Siegel--Weil
formula for similitudes.  In this note, we show that, thanks to the 
recent improvement on the Ramanujan bound
due to Kim--Shahidi \cite{\kimshahidi},
together with a slight variation in the setup of (ii), our method 
yields Jacquet's conjecture in
general.

Since the exposition in \cite{\annals} was specialized from the start 
to the case of
interest for certain arithmetic applications, we will briefly sketch 
the method in general in
the first few sections. We then prove the facts required about the 
extended Siegel--Weil formula.

Several authors have considered interpretations of the vanishing of
the central value $L(\frac12,\pi_1\tt\pi_2\tt\pi_3)$. Here we mention 
only the work of
Dihua Jiang, \cite{\jiang}, who gave an intriguing relation with a 
period of an Eisenstein series on $G_2$.

The authors would like to thank the IHP in Paris where this project 
was realized
during the special program on `geometric aspects of automorphic forms'
in June of 2000.


\subheading{\Sec1. The integral representation of the triple product 
L--function}

Let $G=GSp_6$ be the group of similitudes of the standard $6$ 
dimensional symplectic vector space over $\kay$,
and let $P = MN$ be the Siegel parabolic. For $a\in GL_3$, $b\in 
\Sym_3$ and $\nu$ a scalar, let
$$m(a) = \pmatrix a&{}\\{}&{}^ta^{-1}\endpmatrix, \quad n(b) = 
\pmatrix 1&b\\{}&1\endpmatrix,
\quad\text{and}\quad d(\nu) = \pmatrix 1&{}\\{}&\nu\endpmatrix\in G.\tag1.1$$
Let $K_G = K_{G,\infty}\cdot K_{G,f}$ be the standard maximal compact 
subgroup of $G(\A)$.
For $s\in \C$, let $\lambda_s$ be the character of $P(\A)$ defined by
$$\lambda_s(d(\nu)n(b)m(a)) = |\nu|^{-3s}\,|\det(a)|^{2s}.\tag1.2$$
Let $I(s) = I_P^G(\l_s)$ be the normalized induced representation of 
$G(\A)$, consisting of
all smooth $K_G$--finite functions $\P_s$ on $G(\A)$
such that
$$\P_s(d(\nu)n(b)m(a)g,s) = |\nu|^{-3s-3}\,|\det(a)|^{2s+2}\,\P_s(g).\tag1.3$$
The Eisenstein series associated to a section $\P_s\in I(s)$ is 
defined for $\Re(s)>1$ by
$$E(g,s,\P_s) = \sum_{\gamma\in P(\smallkay)\back G(\smallkay)} 
\P_s(\gamma g),\tag1.4$$
and the normalized Eisenstein series is
$$E^*(g,s,\P_s) = b_G(s)\cdot E(g,s,\P_s),\tag1.5$$
where $b_G(s) = \zeta_\smallkay(2s+2)\,\zeta_{\smallkay}(4s+2),$ as 
in \cite{\psrallis}.
Note that the central character of $E(g,s,\P_s)$ is trivial.
These functions have meromorphic analytic continuations to the whole 
$s$--plane and have no poles on the
unitary axis $\Re(s)=0$. In particular, the map
$$E^*(0):I(0) \lra \Cal A(G),\qquad \P_0 \mapsto (g\mapsto 
E^*(g,0,\P_s))\tag1.6$$
gives a $(\frak g_\infty,K_{G,\infty})\times G(\A_f)$--intertwining map
from the induced representation $I(0)$ at $s=0$ to the space of 
automorphic forms
on $G$ with trivial central character.

Let
$$\align
\bold{G} &= (GL_2\times GL_2\times GL_2)_0\tag1.7\\
\nass
{}& =\big\{(g_1,g_2,g_3)\in (GL_2)^3\mid 
\det(g_1)=\det(g_2)=\det(g_3)\big\}.\endalign
$$
This group embeds diagonally in $G = GSp_6$. For automorphic forms 
$f_i\in \pi_i$, $i=1$, $2$, $3$,
let $F =f_1\tt f_2\tt f_3$ be the corresponding function on $\bold{G}(\A)$.
The global zeta integral \cite{\psrallis} is given by
$$Z(s,F,\P_s) = \int_{Z_{G}(\A)\bold{G}(\smallkay)\back \bold{G}(\A)} 
E^*(\bold{g},s,\P_s)\,F(\bold{g})\,d\bold{g}.\tag1.8$$
Suppose that the automorphic forms $f_i\in \pi_i$ have factorizable 
Whittaker functions
$W_i^\psi = \tt_v W_{i,v}^\psi$ and that the
section $\P_s$ is factorizable.
Let $S$ be a finite set of places of $\kay$, including all 
archimedean places, such that, for $v\notin S$,
\roster
\item"{(i)}" the fixed additive character $\psi$ of $\A/\kay$ has 
conductor $\Cal O_{\smallkay,v}$ at $v$.
\item"{(ii)}" $\pi_{i,v}$ is unramified, $f_i$ is fixed under $K_v = 
GL_2(\Cal O_{\smallkay,v})$,
and $W_{i,v}^\psi(e)=1$.
\item"{(iii)}" $\P_{s,v}$ is right invariant under $G(\Cal 
O_{\smallkay,v})=K_{G,v}$ and $\P_{s,v}(e)=1$,
\endroster
Then
$$Z(s,F,\P_s) = L^S(s+\frac12,\pi_1\tt\pi_2\tt\pi_3)\cdot\prod_{v\in 
S} Z_v(s,W_v^\psi,\P_{s,v}),\tag1.9$$
for local zeta integrals $Z_v(s,W_v^\psi,\P_{s,v})$, where
$W_v^\psi = W_{1,v}^\psi\,W_{2,v}^\psi\,W_{3,v}^\psi$. Here
$$Z(s,W_v^\psi,\P_{s,v}) = \int_{Z_{G}(\smallkay_v)\bold 
M(\smallkay_v)\back \bold{G}(\smallkay_v)}
\P_{s,v}(\delta\, g)\, W_v^\psi(g)\,dg,\tag1.10$$
where $\delta\in G(\kay)$ is a representative for the open orbit of 
$\bold{G}$ in $P\back G$, cf. for example
\cite{\grosskudla}, and
$$\bold M= \{(\pmatrix 1&x_1\\{}&1\endpmatrix,\pmatrix 
1&x_2\\{}&1\endpmatrix,\pmatrix 1&x_3\\{}&1\endpmatrix) \in \bold{G}
\mid x_1+x_2+x_3=0\}.\tag1.11$$
Here $L^S(s,\pi_1\tt\pi_2\tt\pi_3)$ is the triple product L-functions with the
factors for $v\in S$ omitted.

\subheading{\Sec2. Local zeta integrals}

In this section, we record some consequences of recent results of Kim and
Shahidi
\cite{\kimshahidi} on the Ramanujan estimate
for the $\pi_i$'s.  We begin by recalling relevant aspects of the
local theory of the triple product, as recently completed by
Ikeda and Ramakrishnan.  In the following proposition
by ``local Euler factor" at a finite place $v$ of $\kay$ we mean
a function of the form $P(q_v^{-s})^{-1}$, where $P$ is a
polynomial, $P(0) = 1$, and $q_v$ is the order of the residue field; at an
archimedean field we mean a finite product of Tate's local Euler factors
for $GL(1)$.

\proclaim{Proposition 2.1} Let $v$ be a place of $\kay$ and let
$\pi_{i,v}$, $i = 1, 2, 3$, be a triple of admissible irreducible
representations
of $GL(2,\kay_v)$ that arise as local components at $v$ of
cuspidal automorphic representations $\pi_i$.

(i) There exists a local Euler factor
$L(s,\pi_{1,v}\otimes\pi_{2,v}\otimes\pi_{3,v})$ such that, for any
local data $(W_v^\psi,\P_{s,v})$, the quotient
$$\tilde{Z}_v(s,W_v^\psi,\P_{s,v}) = Z_v(s,W_v^\psi,\P_{s,v})\cdot
L(s+\frac12,\pi_{1,v}\tt\pi_{2,v}\tt\pi_{3,v})^{-1}$$
is entire as a function of $s$.

(ii) Let $\sigma_{i,v}$, $i = 1, 2, 3$, be the representations of the
Weil-Deligne group of $\kay_v$ associated to $\pi_{i,v}$ by the
local Langlands correspondence.  Then
$$L(s,\pi_{1,v}\otimes\pi_{2,v}\otimes\pi_{3,v}) =
L(s,\sigma_{1,v}\otimes\sigma_{2,v}\otimes\sigma_{3,v}).$$


(iii) For any finite place $v$, there is a local section
$\P_{s,v}$ and a Whittaker function $W_v^\psi =
W_{1,v}^\psi\,W_{2,v}^\psi\,W_{3,v}^\psi$,  such that
$$Z(s,\P_{s,v},W_v^\psi) \equiv 1.$$

(iv) For any archimedean place $v$, there exists a finite collection of
Whittaker functions $W_v^{\psi,j}$
and of sections $\P_{s,v}^j$, holomorphic in a neighborhood of $s=0$
such that
$$\sum_j Z(0,\P_{s,v}^j,W_v^{\psi,j}) =1.$$
\endproclaim

\demo{Proof}  For $v$ non-archimedean, assertion (i) is proved in \S 3,
Appendix 3, of \cite{\psrallis};
see
\cite{\ikedapoles}, p. 227 for a concise statement.  For $v$ real or
complex, (i) and
(ii) were proved in several steps by Ikeda, of which the crucial
one is \cite{\IkedaDuke}, Theorem 1.10.  Assertion (ii) in
general is due to Ramakrishnan, \cite{\dinakar}, Theorem 4.4.1.
For the moment, the hypothesis that the $\pi_{i,v}$ embed in global
cuspidal representations seems to be necessary.

Assertions (iii) and (iv) are contained in Proposition 3.3 of \cite{\psrallis}.
\qed
\enddemo

\proclaim{Proposition 2.2} (i) For any triple $\pi_i$ of cusp forms for
$GL_2$ over $\kay$, and for any
place $v$, the local Langlands L--factor
$L(s,\pi_{1,v}\tt\pi_{2,v}\tt\pi_{3,v})$ is holomorphic at $s=\frac12$.
\hfill\break
(ii) For any place $v$, for any triple of Whittaker functions
$W_{i,v}^\psi$ in the Whittaker spaces
of $\pi_{i,v}$, and for any section $\P_{s,v}\in I_v(s)$, holomorphic in a
neighborhood of $s=0$,
the local zeta integral
$Z(s,W_v^\psi,P_{s,v})$ is holomorphic in a neighborhood of $s=0$.
\hfill\break
\endproclaim

\demo{Proof}  This follows from the results of Kim and Shahidi.  We sketch
the simple
argument, quoting the proof of Proposition 3.3.2 of \cite{\dinakar}.  Let
$\sigma_{i,v}$ correspond
to $\pi_{i,v}$ as in the previous proposition.  For the present purposes we
can assume
each $\pi_i$ to be unitary.  Indeed, this can be arranged by twisting
$\pi_i$ by a (unique) character
of the form $|\cdot|^{a_i}$, where $|\cdot|$ is the id\`ele norm and
$a_i \in \Bbb{C}$.  Since
the product of the central characters of $\pi_i$ is trivial, we have $a_1 +
a_2 + a_3 = 0$,
so the triple product $L$-factor is left unaffected.

To each $\pi_{i,v}$, necessarily generic and now assumed unitary, we can
assign an index $\lambda_{i,v}$
which measures the failure of $\pi_{i,v}$ to be tempered; we have
$\lambda_{i,v} = t$ if
$\pi_{i,v}$ is a complementary series attached to
$(\mu|\cdot|^{t},\mu|\cdot|^{-t})$ with $t > 0$ and $\mu$ unitary,
$\lambda_{i,v} = 0$ otherwise.  Then, according to \cite{\dinakar}, (3.3.10),
$$L(s, \pi_{1,v}\tt\pi_{2,v}\tt\pi_{3,v})\quad  \text{is holomorphic for}
\Re(s) > \lambda(\pi_{1,v}) + \lambda(\pi_{2,v}) +
\lambda(\pi_{3,v}).  \tag 2.1$$

Now (i) follows from (2.1) and the Kim-Shahidi estimate $\lambda(\pi_{i,v})
< \frac{5}{34}$
for all $i$ and all $v$ \cite{\kimshahidi}, whereas (ii) follows from (i),
and Proposition~2.1 (i) and (ii).
\qed\enddemo

By (1.8), (1.9), and the holomorphy of $E^*(g,s,\P)$ on the unitary
axis,  the expression
$$\align
L^S(s+\frac12,&\pi_1\tt\pi_2\tt\pi_3)\cdot\prod_{v\in S} 
Z_v(s,F,\P_{s,v})\tag2.2\\
\nass
{}&= \int_{Z_{G}(\A)\bold{G}(\smallkay)\back \bold{G}(\A)}
E^*(\bold{g},s,\P_s)\,F(\bold{g})\,d\bold{g}.
\endalign$$
is holomorphic at $s=0$ for all choices of data $F$ and $\P_s$.
By varying the data for places in $S$ and applying (iii) and (iv) of 
Proposition~2.1, it follows
that the partial Euler product $L^S(s+\frac12,\pi_1\tt\pi_2\tt\pi_3)$ 
is holomorphic at $s=0$.
By (i) of Proposition~2.2, the Euler product
$L(s+\frac12,\pi_1\tt\pi_2\tt\pi_3)$ over all finite places is 
holomorphic at $s=0$, and we obtain
the identity
$$\align
L(\frac12,&\pi_1\tt\pi_2\tt\pi_3)\cdot\prod_{v\in S}
Z_v^*(0,W_v^\psi,\P_{s,v}) \tag2.3\\
\nass
{}&=\int_{Z_{G}(\A)\bold{G}(\smallkay)\back \bold{G}(\A)}
E^*(\bold{g},0,\P_s)\,F(\bold{g})\,d\bold{g}.
\endalign$$
where
$$Z_v^*(s,W_v^\psi,\P_{s,v}) =\cases \tilde{Z}_v(s,W_v^\psi,\P_{s,v})
&\text{ if $v\in S_f$,}\\
\nass
Z_v(s,W_v^\psi,\P_{s,v}) &\text{ if $v\in S_\infty$.}
\endcases\tag2.4
$$

\proclaim{Corollary 2.3} $L(\frac12,\pi_1\tt\pi_2\tt\pi_3) = 0$ if and only if
$$ \int_{Z_{G}(\A)\bold{G}(\smallkay)\back \bold{G}(\A)}
E^*(\bold{g},0,\P_s)\,F(\bold{g})\,d\bold{g}=0,$$
for all choices of $F\in \Pi = \pi_1\tt\pi_2\tt\pi_3$
and $\P_s\in I(s)$.
\endproclaim

Of course, relation (2.3) gives a formula for
$L(\frac12,\pi_1\tt\pi_2\tt\pi_3)$ for a
suitable choice of $F$ and $\P_s$.

\subheading{\Sec3. The Weil representation for similitudes}

The material of this section is a slight variation on that of section 5 of
\cite{\duke}. We consider only the case of the dual pair 
$(GO(V),GSp_6)$ where the
space $V$ has square discriminant.

Let $B$ be a quaternion algebra over $\kay$ (including the 
possibility $B=M_2(\kay)$), and let
Let $V=B$ be a $4$ dimensional quadratic space over $\kay$ where the 
quadratic form
is given by $Q(x) = \a\, \nu(x)$, where $\nu$ is the reduced norm on 
$B$ and $\a\in \kay^\times$.
Note that the isomorphism class of $V$ is determined by $B$ and the 
sign of $\a$ at the
set $\Sigma_\infty(B)$ of real archimedean places of $\kay$ at which 
$B$ is division.
Let $H=GO(V)$
and let $H_1= O(V)$ be the kernel of the scale map $\nu:H\rightarrow \Bbb G_m$.
Let $G = GSp_6$, and let $G_1=Sp_6$ be the kernel of the scale map 
$\nu:G\rightarrow \Bbb G_m$.
We want to extend the standard Weil representation $\o = \o_\psi$ of 
$H_1(\A)\times G_1(\A)$
on the Schwartz space $S(V(\A)^3)$. First,
there is a natural action of $H(\A)$ on $S(V(\A)^3)$ given by
$$L(h)\ph(x) = |\nu(h)|^{-3} \,\ph(h^{-1}x).\tag3.1$$
For $g_1\in G_1(\A)$ one has
$$L(h)\o(g_1)L(h)^{-1} = \o(d(\nu)g_1 d(\nu)^{-1}),\tag3.2$$
where $\nu =\nu(h)$, and $d(\nu)$ is as in section 1.
Therefore, one obtains a representation of the semidirect product
$H(\A)\ltimes G_1(\A)$ on $S(V(\A)^3)$. Let
$$R = \{ (h,g)\in H\times G\mid \nu(h)=\nu(g)\}.\tag3.3$$
Then there is an isomorphism
$$R \lra H\ltimes G_1, \qquad (h,g) \mapsto (h,d(\nu(g))^{-1}g) = 
(h,g_1),\tag3.4$$
(this defines a map $g\mapsto g_1$) and a representation of $R(\A)$ 
on $S(V(\A)^3)$
given by
$$\o(h,g)\ph(x) = (L(h)\o(g_1)\ph)(x) = 
|\nu(h)|^{-3}(\o(g_1)\ph)(h^{-1}x).\tag3.5$$
The theta distribution $\Theta$ on $S(V(\A)^3)$ is invariant under 
$R(\kay)$, since,
for $(h,g)\in R(\kay)$,
$$\align
\Theta(\o(h,g)\ph) & = \sum_{x\in V(\smallkay)^3} 
|\nu(h)|^{-3}(\o(g_1)\ph)(h^{-1}x)\\
\nass
{}&= \sum_{x\in V(\smallkay)^3} (\o(g_1)\ph)(x)\tag3.6\\
\nass
{}&= \Theta(\o(g_1)\ph) = \Theta(\ph),
\endalign
$$
since $g_1\in G_1(\kay)$.
The theta kernel, defined for $(h,g)\in R(\A)$ by
$$\theta(h,g;\ph) = \sum_{x\in V(\smallkay)^3} \o(h,g)\ph(x),\tag3.7$$
is thus left $R(\kay)$ invariant.

{\bf Remark 3.1.} Aside from a shift in notation, the convention here is
essentially the same as in section 5 of \cite{\duke} and section 3 of 
\cite{\annals},
{\it except that} we take pairs $(h,g)$ here versus $(g,h)$ there. 
Compare (3.5)
above with (5.1.5) of \cite{\duke}. It turns out that this seemingly 
slight shift in convention will
crucial for the extension of the Siegel--Weil formula to similitudes, 
as we will
see below.

Note that the set of archimedean places $\Sigma_\infty(B)$ introduced above
is the set of all real archimedean places of $\kay$ at which
$V$ is definite (positive or negative). Then,
$$\align
G(\A)^+ &:= \{g\in G(\A)\mid \nu(g) \in \nu(H(\A))\}\tag3.8\\
\nass
{}&=\{g\in G(\A)\mid \nu(g)_v>0, \ \forall v\in \Sigma_\infty(V)\}.\endalign$$
For $g\in G(\A)^+$, and $\ph\in S(V(\A)^3)$, and for $V$ anisotropic 
over $\kay$,
i.e., for $B\ne M_2(\kay)$, the theta integral is defined by
$$I(g,\ph) = \int_{H_1(\smallkay)\back H_1(\A)} 
\theta(h_1h,g;\ph)\,dh_1,\tag3.9$$
where $h\in H(\A)$ with $\nu(h)=\nu(g)$. It does not depend on the 
choice of $h$.

In the case $B=M_2(\kay)$, the theta integral must be defined by 
regularization.
If $\kay$ has a real place, the procedure outlined on p.621 of \cite{\annals},
\cite{\krannals}, using a certain differential operator to kill 
support, can be applied.
An analogous procedure using
an element of the Bernstein center can be applied at a nonarchimedean place,
\cite{\vtan}. We omit the details.

\proclaim{Lemma 3.2} (i) (Eichler's norm Theorem) If $\a\in 
\nu(H(\A))\cap \kay^\times$,
then there exists an element $h\in H(\kay)$ with $\nu(h) = \a$. \hfill\break
(ii) The theta integral is left invariant under $G(\A)^+\cap 
G(\kay)$. \hfill\break
(iii) The theta integral has trivial central character, i.e.,
for $z\in Z_G(\A)\subset G(\A)^+$, $I(zg,\ph) = I(g,\ph)$.
\endproclaim
\demo{Proof} (i) is a standard characterization of $\nu(H(\kay))$ is 
the present case.
To check (ii), given $\gamma\in G(\A)^+\cap G(\kay)$, choose 
$\gamma'\in H(\kay)$
with $\nu(\gamma')=\nu(\gamma)$. Then
$$\align
I(\gamma g, \ph)&= \int_{H_1(\smallkay)\back H_1(\A)} 
\theta(h_1\gamma'h,\gamma g;\ph)\,dh_1\\
\nass
{}&= \int_{H_1(\smallkay)\back H_1(\A)} \theta(\gamma'h_1h,\gamma 
g;\ph)\,dh_1\tag3.10\\
\nass
{}& = I(g,\ph),
\endalign
$$
via the left invariance of the theta kernel under 
$(\gamma',\gamma)\in R(\kay)$. Here, in the next to last
step, we have conjugated the domain of integration $H_1(\kay)\back 
H_1(\A)$ by the element
$\gamma'\in H(\kay)$.

Finally, the proof of (iii) is just like that of Lemma~5.1.9 (ii) in 
\cite{\duke}.
\qed\enddemo

Since $G(\A) = G(\kay)G(\A)^+$, it follows that $I(g,\ph)$ has a 
unique extension to a left
$G(\kay)$--invariant function on $G(\A)$. Moreover, for any $g_0\in 
G(\A)^+$, we have
$$I(gg_0,\ph) = I(g, \o(h_0,g_0)\ph),\tag3.11$$
where $h_0\in H(\A)$ with $\nu(h_0)=\nu(g_0)$. In particular, if 
$h_1\in H_1(\A)$, then
$$I(g,\o(h_1)\ph) = I(g,\ph).\tag3.12$$

\subheading{\Sec4. The Siegel--Weil formula for $(GO(V),GSp_6)$ }

First we recall the Siegel--Weil formula for $(O(V),Sp_6)$. The 
results of \cite{\krannals}
on the regularized Siegel--Weil formula were formulated over a 
totally real number field, since,
at a number of points, we needed facts about degenerate principal 
series, intertwining
operators, etc. which had not been checked for complex places. The 
proof in the case
of the central value of the Siegel--Eisenstein series is simpler than 
the general case,
and the additional facts needed at complex places are easy to check. In the
rest of this section, we will state the results for an arbitrary 
number field $\kay$.
A sketch of the proof of Theorem~4.1 below for such a field $\kay$ 
will be given in the Appendix below.

Let $I_1(s) = I_{P_1}^{G_1}(\l_s)$ be the global induced representation of
$G_1(\A) = Sp_6(\A)$ induced from the restriction of the character $\l_s$
of $P(\A)$ to $P_1(\A) = P(\A)\cap G_1(\A)$.

For a global quadratic space $V$ of dimension $4$ over $\kay$ associated to a
quaternion algebra $B$, as in the previous section, there is a
$(\frak g_{1,\infty},K_{G_1,\infty})\times G_1(\A_f)$--equivariant map
$$S(V(\A)^3) \lra I_1(0),\qquad \ph\mapsto [\ph],\tag4.1$$
where
$$[\ph](g_1) = (\o(g_1)\ph)(0).\tag4.2$$
The image, $\Pi_1(V)$, is an irreducible summand of the unitarizable induced
representation $I_1(0)$. By the results of Rallis \cite{\rallisHDC}, 
Kudla--Rallis \cite{\krrdps}, \cite{\krinvardist},
and the appendix,
$$\Pi_1(V)\simeq  S(V(\A)^3)_{O(V)(\A)},\tag4.3$$
the space of $H_1(\A)=O(V)(\A)$--coinvariants.
One then has a decomposition
$$I_1(0) = \bigg(\oplus_{V} \Pi_1(V)\bigg)\oplus\bigg(\oplus_{\Cal V} 
\Pi_1(\Cal V)\bigg),\tag4.4$$
into irreducible representations of $G_1(\A)$, as $V$ runs over the 
isomorphism classes of
such spaces and as $\Cal V$ runs over the
incoherent collections, obtained by switching one local component of 
a $\Pi_1(V)$, cf. \cite{\annalsII}.

The Siegel--Weil formula of \cite{\annals}, asserts the following in 
the present case.

\proclaim{Theorem 4.1} (i) The $(\frak 
g_{1,\infty},K_{G_1,\infty})\times G_1(\A_f)$--intertwining map
$$E_1(0):I_1(0)\lra \Cal A(G_1), \qquad \P_0\mapsto (g_1\mapsto 
E(g_1,0,\P_s))$$
has kernel $\oplus_{\Cal V} \Pi_1(\Cal V)$. \hfill\break
(ii) For a section $\P_s\in I_1(s)$ with $\P_0=[\ph]$ for some 
$\ph\in S(V(\A)^3)$,
$$E(g_1,0,\P_s) = 2 \,I(g_1,\ph),\tag SW$$
for the theta integral as defined in \Sec3.
\endproclaim

As explained in the previous section, the theta integral can be extended to
an automorphic form on $G(\A)$. We will see presently that it 
coincides with  an Eisenstein series on $G(\A)$.

Restriction of functions
from $G(\A)$ to $G_1(\A)$ yields an isomorphism $I(s) \isoarrow I_1(s)$, which
is intertwining for the right action of $G_1(\A)$.  Here $I(s)$ is the
induced representation of $G(\A)$ defined in section 1.
The inverse map  is given by $\P_s\mapsto \P_s^\sim$ where
$$\P_s^\sim(g) = |\nu(g)|^{-3s-3}\P_s(g_1),\tag4.5$$
for $g_1=d(\nu(g))^{-1}g$, as above. The decomposition (4.4) into 
$G_1(\A)$--irreducibles yields a decomposition
$$I(0) = \bigg(\oplus_{B} \Pi(B)\bigg)\oplus\bigg(\oplus_{\Cal B} 
\Pi(\Cal B)\bigg),\tag4.6$$
into irreducible representations of $G(\A)$, where, for a global 
quaternion algebra $B$
over $\kay$,
$$\Pi(B) = \oplus_V \Pi(V)\tag4.7$$
where $V$ runs over the non-isomorphic spaces associated to $B$ 
(i.e., different multiples
of the norm form) and $\Pi(V)$ denotes the image of $\Pi_1(V)$
under the inverse of the restriction isomorphism. Note that there are 
$2^{|\Sigma_\infty(B)|}$ such $V$'s.
In effect, at a real archimedean place $v$, the local induced 
representation $I_1(0)_v$ has
a decomposition into irreducible $(\frak g_{1,v},K_{G_1,v})$--modules
$$I(0)_v = \Pi(4,0)_v\oplus\Pi(2,2)_v\oplus\Pi(0,4)_v\tag4.8$$
according to signatures. The space $\Pi(2,2)_v$ is actually stable 
under $(\frak g_v,K_{G,v})$, as is
the sum $\Pi(4,0)_v\oplus\Pi(0,4)_v$, and
$$\Pi(B)_v = \cases \Pi(2,2)_v &\text{ if $B_v\simeq M_2(\R)$.}\\
\nass
\Pi(4,0)_v\oplus\Pi(0,4)_v&\text{ if $B_v$ is division. }
\endcases\tag4.9
$$
The summands $\Pi(\Cal B)$ are defined similarly.

\proclaim{Theorem 4.2} (i) The $(\frak g_{\infty},K_{G,\infty})\times 
G(\A_f)$--intertwining map
$$E(0):I(0)\lra \Cal A(G), \qquad \P_s\mapsto (g\mapsto E(g,0,\P_s))$$
has kernel $\oplus_{\Cal B} \Pi(\Cal B)$. \hfill\break
(ii) For a section $\P_s\in I(s)$ with $\P_0\in \Pi(V)$ so that 
$\P_0=[\ph]^\sim$ for some $\ph\in S(V(\A)^3)$,
$$E(g,0,\P_s) = 2 \,I(g,\ph),\tag GSW$$
for the theta integral as defined in \Sec3.
\endproclaim

\demo{Proof} For $g_0\in G(\A_f)$, we have
$$E(gg_0,s,\P_s) = E(g,s,r_s(g_0)\P_s),\tag4.10$$
where $r_s$ denotes the action in the induced representation $I(s)$ 
by right translation.
Taking the value at $s=0$, we obtain
$$E(gg_0,0,\P_s) = E(g,0,r_s(g_0)\P_s).\tag4.11$$
Note that this value depends only on $\P_0$ and $r_0(g_0)\P_0$.

\proclaim{Lemma 4.3} For $\ph\in S(V(\A)^3)$, let $[\ph]\in I_1(0)$
be defined by (4.2) and let $[\ph]^\sim$ be the corresponding 
function in $I(0)$ under the
inverse of the restriction isomorphism.
\hfill\break
(i) For $g\in G(\A)^+$,
$$[\ph]^\sim(g) = (\o(h,g)\ph)(0),$$
where $h\in GO(V)(\A)$ with $\nu(h)=\nu(g)$.
\hfill\break
(ii) For $g_0\in G(\A_f)$,
$$r_0(g_0)[\ph]^\sim = [\o(h_0,g_0)\ph]^\sim,$$
where $h_0\in GO(V)(\A_f)$ with $\nu(h_0)=\nu(g_0)$.
\endproclaim

\demo{Proof} For (i), we have
$$\align
[\ph]^\sim(g) &= [\ph]^\sim(d(\nu) g_1)\\
\nass
{}&= |\nu|^{-3} [\ph](g_1)\\
\nass
{}&=|\nu|^{-3}(\o(g_1)\ph)(0)\tag4.12\\
\nass
{}&= (L(h)\o(g_1)\ph)(0)\\
\nass
{}&=(\o(h,g)\ph)(0).
\endalign
$$
For (ii),
$$\align
(r_0(g_0)[\ph]^\sim)(g) &= |\nu|^{-3}[\ph]^\sim(g_1g_0)\\
\nass
{}&=|\nu|^{-3}(\o(h_0,g_1g_0)\ph)(0)\\
\nass
{}&=|\nu|^{-3}(\o(g_1)\o(h_0,g_0)\ph)(0)\tag4.13\\
\nass
{}&= |\nu|^{-3}[\o(h_0,g_0)\ph](g_1)\\
\nass
{}&= [\o(h_0,g_0)\ph]^\sim(g).
\endalign
$$

\qed\enddemo

Thus, if $\P_0 = [\ph]^\sim$, then $r_0(g_0)\P_0 = [\o(h_0,g_0)\ph]^\sim$.
Since $G(\A) = G(\kay)Z_G(\A) G_1(\A) G(\A_f)$, we have, by (4.11),
$$\align
E(g,0,\P_s) &= E(\gamma z g_1g_0,0,\P_s)\tag4.14\\
\nass
{}&= E(g_1,0,r_s(g_0)\P_s),
\endalign
$$
the value at $g_1$ of the Siegel--Eisenstein series attached to 
$\o(h_0,g_0)\ph\in S(V(\A)^3)$.
On the other hand, by (3.11),
$$\align
I(g,\ph)&=I(\gamma z g_1g_0,\ph)\tag4.15\\
\nass
{}&=I(g_1,\o(h_0,g_0)\ph).
\endalign
$$
Thus (GSW) follows from (SW).
\qed\qed\enddemo

\vfill\eject

\subheading{\Sec5. Proof of Jacquet's conjecture}

Applying the Siegel--Weil formula for similitudes to the basic 
identity (2.3), we obtain
$$\align
&L(\frac12,\pi_1\tt\pi_2\tt\pi_3)\cdot Z^*(F,\P)\\
\nass
{}&=\int_{Z_{\bold{G}}(\A)\bold{G}(\smallkay)\back \bold{G}(\A)} 
E^*(\bold{g},0,\P_s)\,F(\bold{g})\,d\bold{g}\tag5.1\\
\nass
{}&=2\zeta_{\smallkay}(2)^2
\sum_{V}\int_{Z_{\bold{G}}(\A)\bold{G}(\smallkay)\back 
\bold{G}(\A)}I(\bold{g},\ph^V)\,F(\bold{g})\,d\bold{g}.
\endalign
$$
where
$$Z^*(F,\P) = \prod_{v\in S} Z_v^*(0,W_v^\psi,\P_{s,v}),\tag5.2$$
and
where $\ph^V\in S(V(\A)^3)$, and, in fact, only a finite set of $V$'s 
occurs in the sum.
More precisely,
in the decomposition (4.6),
$$\P_0 = \sum_V\ [\ph^V]^\sim + \text{ terms in the $\Pi(\Cal B)$'s } 
\in I(0),\tag5.3$$
where the quaternion algebras $B$ associated to $V$'s are split 
outside the set $S$, due
to condition (iii) in the definition of $S$ in section 1.

Now consider the integral in the last line of (5.1) for a fixed $\ph= 
\ph^V$. To apply the
seesaw identity, we set
$$\align
H&=GO(V)\\
\nass
\bold{H}&=\{(h_1,h_2,h_3)\in H^3\mid \nu(h_1)=\nu(h_2)=\nu(h_3)\},\tag5.4\\
\nass
\bold{R}&=\{(\bold{h},\bold{g})\in \bold{H}\times\bold{G}\mid 
\nu(\bold{h})=\nu(\bold{g})\}\\
\nass
R_0&=\{(h,\bold{g})\in H\times \bold{G}\mid \nu(h)=\nu(\bold{g})\},
\endalign
$$
and hence have the seesaw pair:
$$\matrix
I(\cdot,\ph;F)\quad&(GO(V)^3)_0=&\bold{H}&{}&{}&G&=GSp_6&\quad I(\cdot,\ph)\\
\nass
{}&{}&\uparrow&\nwarrow&\nearrow&\uparrow&{}&{}\\
\noalign{\vskip -3pt}
{}&{}&\vert&\swarrow&\searrow&\vert&{}&{}\\
\nass
\triv\quad&GO(V)=&H&{}&{}&\bold{G}&=(GL_2^3)_0&\quad F\endmatrix\tag5.5
$$
There are representations of both $R(\A)$ and $\bold{R}(\A)$ on 
$S(V(\A)^3)$, and the restriction of these
representations to the common subgroup $R_0(\A)$ coincide.

For $F$ a cuspidal automorphic form on $\bold{G}(\A)$ and for 
$\bold{h}\in \bold{H}(\A)$, let
$$I(\bold{h},\ph;F) = \int_{\bold{G}_1(\smallkay)\back \bold{G}_1(\A)}
\theta(\bold{h},\bold{g}_1\bold{g};\ph)\,F(\bold{g}_1\bold{g})\,d\bold{g}_1,\tag5.6$$
where $\bold{g}\in\bold{G}(\A)$ with $\nu(\bold{g})=\nu(\bold{h})$.
\proclaim{Lemma 5.1}{\rm (Seesaw identity)}
$$\int_{Z_{G}(\A)\bold{G}(\smallkay)\back \bold{G}(\A)} 
I(\bold{g},\ph)\,F(\bold{g})\,d\bold{g}
= \int_{Z_H(\A)H(\smallkay)\back H(\A)} I(h,\ph;F)\,dh.
$$
\endproclaim
\demo{Proof}
Note that
$
Z_{G}(\A)\bold{G}(\kay)\back \bold{G}(\A) \simeq 
Z_{G}(\A)\bold{G}(\kay)^+\back \bold{G}(\A)^+,
$
and that
$$Z_{G}(\A)\bold{G}(\kay)^+\bold{G}_1(\A)\back \bold{G}(\A)^+
\simeq Z_H(\A)H(\kay)\back H(\A) \simeq 
\A^{\times,2}\kay^{\times,+}\back \A^{\times,+} =:C, \tag5.7
$$
is compact, where $\A^{\times,+} = \nu(H(\A))$ and $\kay^{\times,+} = 
\nu(H(\kay))$.
Fixing a Haar measure $dc$ giving $C$ volume $1$, we have
$$\align
&\int_{Z_{\bold{G}}(\A)\bold{G}(\smallkay)\back 
\bold{G}(\A)}I(\bold{g},\ph)\,F(\bold{g})\,d\bold{g}\\
\nass
{}&=\int_C\int_{\bold{G}_1(\smallkay)\back \bold{G}_1(\A)}
\int_{H_1(\smallkay)\back H_1(\A)} 
\theta(h_1h(c),\bold{g}_1\bold{g}(c);\ph)\,F(\bold{g}_1\bold{g}(c))\,dh_1\,d\bold{g}_1\,dc\tag5.8\\
\nass
{}&= \int_{Z_H(\A)H(\smallkay)\back H(\A)} I(h,\ph;F)\,dh,
\endalign
$$
generalizing the proof of Proposition~7.1.4 of \cite{\duke}.
\qed\enddemo

To apply the seesaw identity to the restriction to $\bold{G}(\A)$ of 
a function $F\in \Pi = \pi_1\tt\pi_2\tt\pi_3$,
we recall the description, from sections 7 and 8 of \cite{\annals}, 
of the corresponding space of functions $\Theta(\Pi)$
on $\bold H(\A)$ spanned by the $I(\bold h,\ph;F)$'s for $F\in \Pi$ and
$\ph\in S(V(\A)^3)$. Note that one obtains the same space by fixing a 
nonzero $F$ and only varying $\ph$
\cite{\howeps}.

The action of $B^\times\times B^\times$ on $V=B$,
$\rho(b_1,b_2)x = b_1 x b_2^{-1}$
determines an extension
$$ 1\lra \Bbb G_m \lra \big(\, B^\times\times 
B^\times\,\big)\rtimes\langle\bold t\rangle  \lra 
H=GO(V)\lra1\tag5.9$$
where the involution $\bold t$ acts on $V$ by $\rho(\bold t)(x) = x^\iota$
and on $B^\times\times B^\times$ by $(b_1,b_2)\mapsto 
(b_2^\iota,b_1^\iota)^{-1}$.
Write
$$\tilde H = \big(\, B^\times\times B^\times\,\big)\rtimes\langle\bold t\rangle
\qquad\quad\text{and}\qquad\quad \tilde H^0 = B^\times\times 
B^\times,\tag5.10$$
and let $\tilde\bold H$ and $\tilde\bold H^0$ be the analogous groups 
for $\bold H = (GO(V)^3)_0$.
Thus, we have the diagram
$$ \matrix
{}&{}& \tilde\Theta(\Pi)&{}&\Theta(\Pi)\\
\nass
\nass
\tilde\bold H^0 & \hookrightarrow & \tilde \bold H& \lra & \bold H\\
\nass
\uparrow&{}&\uparrow&{}&\uparrow\\
\nass
\tilde H^0&\hookrightarrow& \tilde H & \lra& H
\endmatrix\tag5.11
$$
For an irreducible cuspidal automorphic representation $\pi$ of $GL_2(\A)$,
let $\pi^B$ be the associated automorphic representation of 
$B^\times(\A)$ under
the Jacquet-Langlands correspondence. We take $\pi^B$ to be zero if 
$\pi$ does not
correspond to a representation of $B^\times(\A)$.  Similarly, let
$\Pi^B = \pi_1^B\tt\pi_2^B\tt\pi_3^B$ be the corresponding representation of
$B^\times(\A)^3$, or zero if some factor does not exist. Note that the central
character of $\Pi^B$ is trivial, and so,
$(\Pi^B)^\vee \simeq \Pi^B$, where $(\Pi^B)^\vee$ is the
contragradient of $\Pi^B$. Thus we can view the space of of functions 
$\Pi^B$ on $B^\times(\A)^3$
as the automorphic realization of both $\Pi^B$ and its contragredient.

The following result is proved in \cite{\annals}, sections 7 and 8, 
based on the work of Shimizu and
Prasad.

\proclaim{Proposition 5.2} (i) $\Theta(\Pi)$ is either zero or a 
cuspidal automorphic representation of
$\bold H(\A)$ and is nonzero if and only if $\Pi^B$ is nonzero. \hfill\break
(ii) As spaces of functions on $\tilde\bold H^0(\A)$,
$$ \tilde\Theta(\Pi)\big\vert_{\tilde\bold H^0(\A)} = \bigg(\Pi^B \tt 
(\Pi^B)^\vee\,\bigg)\bigg\vert_{\tilde\bold H^0(\A)}.$$
\endproclaim

For fixed $F\in \Pi$ and $\ph\in S(V(\A)^3)$, we let
$\tilde I(\cdot,\ph;F)$ denote the pullback of $I(\cdot ,\ph;F)$ to 
$\tilde\bold H(\A)$,
and, via (ii) of Proposition~5.2, we write the restriction of this 
function to $\tilde \bold H^0(A)$ as
$$\tilde I((\bold b_1,\bold b_2),\ph;F) = \sum_r I^{1,r}(\bold 
b_1,\ph;F) \,I^{2,r}(\bold b_2,\ph;F)\tag5.12$$
for functions $I^{i,r}(\cdot,\ph;F)\in \Pi^B$ and $\bold b_i\in 
B^\times(\A)^3$.
The seesaw then gives
$$\align
&\int_{Z_{G}(\A)\bold{G}(\smallkay)\back \bold{G}(\A)} 
I(\bold{g},\ph)\,F(\bold{g})\,d\bold{g}\\
\nass
{}&= \int_{Z_H(\A)H(\smallkay)\back H(\A)} I(h,\ph;F)\,dh\\
\nass
{}&=\int_{Z_{\tilde H^0}(\A)\tilde H(\smallkay)\back \tilde H(\A)} 
\tilde I(h,\ph;F)\,dh\tag5.13\\
\nass
{}&=\int_{Z_{\tilde H^0}(\A)\tilde H^0(\smallkay)\back \tilde 
H^0(\A)} \tilde I(h,\ph;F)\,dh\\
\nass
{}&=\sum_r\int_{\A^\times B^\times(\smallkay)\back B^\times(\A)} 
I^{1,r}(b_1,\ph;F)\,db_1
\cdot \int_{\A^\times B^\times(\smallkay)\back B^\times(\A)} 
I^{2,r}(b_2,\ph;F)\,db_2
\endalign
$$
The fact that the integral over $Z_{\tilde H^0}(\A)\tilde 
H(\smallkay)\back \tilde H(\A)$ in the third line
can be replaced by the integral over $Z_{\tilde H^0}(\A)\tilde 
H^0(\smallkay)\back \tilde H^0(\A)$
is (7.3.2), p.632 of \cite{\annals}; its proof in section 8.6, p.636 
of \cite{\annals} depends
on Prasad's uniqueness theorem for invariant trilinear forms.

Finally, we observe that the integrals in this last expression are 
finite linear combinations
of the quantities $I(f_1^B,f_2^B,f_3^B)$ of (0.1) and that, by (ii) 
of Proposition~5.2 ,
every such quantity can be obtained as an integral
$\int_{\A^\times B^\times(\smallkay)\back B^\times(\A)} 
I^{1,r}(b_1,\ph;F)\,db_1$
for some $\ph$, $F$ and $r$.

This finishes the proof of Jacquet's conjecture.

{\bf Remark 5.3.} In fact, by Prasad's uniqueness theorem, if the root number
$$\e(\frac12,\pi_1\tt\pi_2\tt\pi_3) = 1,\tag5.14$$
then there is a unique $B$ for which $\Pi^B\ne0$ and for which
the space of global invariant trilinear forms on $\Pi^B$ has dimension $1$.
The {\it automorphic} trilinear form is given by integration over 
$\A^\times B^\times(\kay)\back B^\times(\A)$
is then non-zero if and only if  $L(\frac12,\pi_1\tt\pi_2\tt\pi_3)\ne0$.
Choose $f_i^B\in \pi_i^B$, $i=1$, $2$, $3$, such that
$$I(f_1^B,f_2^B,f_3^B)\ne0.\tag5.15$$
For any nonzero $F\in \Pi$, we can choose $\ph\in S(V(\A))^3$ such that
$$\tilde I((\bold b, \bold b'),\ph;F) = 
f_1^B(b_1)f_2^B(b_2)f_3^B(b_3) 
\,f_1^B(b'_1)f_2^B(b'_2)f_3^B(b'_3),\tag5.16$$
where $\bold b = (b_1,b_2,b_3)$ and $\bold b'= (b'_1,b'_2,b'_3)$.
We then obtain
$$
L(\frac12,\pi_1\tt\pi_2\tt\pi_3)\cdot Z^*(F,\P) = 
2\zeta_{\smallkay}(2)^2\,I(f_1^B,f_2^B,f_3^B)^2.\tag5.17
$$
where, $\P$ is determined by $\ph$, and $Z^*(F,\P) \ne0$. Of course, 
this identity
is only useful when one has sufficient information about the function 
$\ph$ and the product
of local zeta integrals $Z^*(F,\P)$.  This was a main concern in 
\cite{\annals}.

On the other hand, when the root number 
$\e(\frac12,\pi_1\tt\pi_2\tt\pi_3) = -1$, then
there is no $\Pi^B$ which supports an invariant trilinear form, and
the central value of the triple product L-function vanishes due to 
the sign in the functional equation.

\subheading{ Appendix: The Siegel--Weil formula for general $\kay$}

In this appendix, we will sketch the proof of Theorem~4.1 for an 
arbitrary number field
$\kay$, indicating the additional facts which are needed when $\kay$ 
has complex places.

First, suppose that $v$ is a complex place of $\kay$ and
consider the local degenerate principal series representation
$I_{1,v}(0)$ of $G_{1,v} = Sp_3(\C)$ and the Weil representation of 
$G_{1,v}$ on $S(V_v^3)$, where
$V_v \simeq M_2(\C)$ with $Q(x)=\det(x)$.
\proclaim{Lemma A.1} (i) $I_{1,v}(0)$ is an irreducible unitarizable 
representation of $G_{1,v}$.
\hfill\break
(ii) (Coinvariants) The map
$S(V_v^3) \rightarrow I_{1,v}(0)$, $\ph\mapsto [\ph]$, analogous to 
(4.1) induces
an isomorphism
$$S(V_v^3)_{H_{1,v}} \isoarrow I_{1,v}(0).$$
Here $H_1=O(V)$.
\endproclaim

{\bf Remark:} Statement (i) is in
Sahi's paper, \cite{\sahi}, Theorem 3A. The proof of (ii)
was explained to us by Chen-bo Zhu\footnote{He also directed us to 
\cite{\sahi}.
We wish to thank him for his help on these points.}  \cite{\zhuII}, 
and is based on
the method of \cite{\zhuI}.

In the case $B=M_2(\kay)$, the theta integral must be defined by 
regularization,
-- cf. the remarks before Lemma~3.2 above.
can be applied. Alternatively,
Siegel--Weil formula for certain
We write
$$I_{\roman{reg}}(g_1,\ph) = \cases I(g_1,\ph)&\text{ if $V$ is anisotropic,}\\
\nass
B_{-1}(g_1,\ph)&\text{ if $V$ is isotropic},
\endcases\tag A.1$$
where $B_{-1}$ is as in (5.5.24) of \cite{\krannals}, except that
we normalize the auxillary Eisenstein series $E(h,s)$ to have
residue $1$ at $s_0'$. The key facts which we need are
the following.
\proclaim{Lemma A.2} (i) The map $I_{\roman{reg}}:S(V(\A)^3) 
\rightarrow \Cal A(G_1)$
factors through the space of coinvariants $S(V(\A)^3)_{H_1(\A)} = 
\Pi(V)$. \hfill\break
(ii) For all $\beta\in \roman{Sym}_3(\kay)$,
$\beta\in \roman{Sym}_3(\kay)$
$$I_{\roman{reg},\beta}(g_1,\ph) = \frac12\cdot \int_{H_1(\A)} 
\o(g_1)\ph(h^{-1}x)\,dh,$$
where $x\in V(\kay)^3$ with $Q(x) =\beta$.
\endproclaim
The second statement here is Corollary~6.11 of \cite{\krannals}; it 
asserts that
the nonsingular Fourier coefficients behave as though no regularization
were involved.

Next we have the analogue of Lemma~4.2, p.111 of \cite{\rallisbook}; the main
point of the proof is the local uniqueness, and the `submersive set' argument
for the archimedean places carries over for a complex place.

\proclaim{Lemma A.3} For $\beta\in \roman{Sym}_3(\kay)$ with
$\det(\beta)\ne0$, let $\Cal T_\beta$ be the space of
distributions $T\in S(V(\A)^3)'$ such that
\hfill\break
(i) $T$ is $H(\A)$--invariant.
\hfill\break
(ii) For all $b\in \roman{Sym}_3(\A_f)$,
$$T(\o(n(b)\ph) = \psi_\beta(b)\,T(\ph),$$
where $\psi_\beta(b) = \psi(\tr(\beta b)).$
\hfill\break
(iii) For an archimedean place $v$ of $\kay$ and for
all $X\in \frak n=\roman{Lie}(N)$,
$$T(\o(X)\ph) = d\psi_\beta(X)\cdot T(\ph).$$
Then $\Cal T_\beta$ has dimension at most $1$ and is spanned by
the orbital integral
$$T(\ph) = \int_{H_1(\A)} \ph(h^{-1}x)\,dh,$$
where $x\in V(\kay)^3$ with $Q(x) = \beta$. In particular, $\Cal T_\beta=0$
if an only if there is no such $x$.
\endproclaim

\demo{sketch of the Proof of Theorem~4.1}
First consider a global space $V$ associated to a quaternion algebra 
$B$. We have
two intertwining maps
$$E_1(0): \Pi(V) \lra \Cal A(G_1) \qquad\text{\rm and }\qquad 
I_{\roman{reg}}:\Pi(V) \lra \Cal A(G_1)\tag A.2$$
from the irreducible representation $\Pi(V) \simeq 
S(V(\A)^3)_{H(\A)}$ of $G_1(\A)$ to the space of
automorphic forms.  For a nonsingular $\beta\in \roman{Sym}_3(\kay)$, 
the distributions
obtained by taking the $\beta$th Fourier coefficient of the 
composition of the projection $S(V(\A)^3)\rightarrow \Pi(V)$
with each of the embeddings in (A.2) satisfy the conditions of 
Lemma~A.3 and hence are proportional. In particular, the
$\beta$-th Fourier of the Eisenstein series vanishes unless $\beta$ 
is represented by $V$.
By the argument of
pp. 111--115 of \cite{\rallisbook}, the constant of proportionality 
is independent of $\beta$
and so there is a constant $c$ such that $E_1(g,0,[\ph]) - c\cdot 
I_{\roman{reg}}(g,\ph)$
has vanishing nonsingular Fourier coefficients.  But then the 
argument at the top of p.28 of \cite{\krannals},
cf. also, \cite{\rallisbook}, implies that this difference is identically zero.

In the case of a component $\Pi(\Cal V)\subset I_1(0)$, the nonsingular Fourier
coefficients of $E(g,0,\P)$ vanish by the argument on p. 28 of 
\cite{\krannals}, so,
again by `nonsingularity' the map $E_1(0)$ must vanish on $\Pi(\Cal V)$.
\qed
\enddemo

\redefine\vol{\oldvol}

\medskip
\centerline{Michael Harris \hfill Stephen S. Kudla}
\centerline{Institut de Math\'ematiques de Jussieu, U.M.R. 7586 du CNRS\hfill Department of Mathematics}
\centerline{175 rue du Chevaleret\hfill University of Maryland}
\centerline{Paris 75013 \hfill College Park, Maryland, 20817}
\centerline{France \hfill U.S.A}
\centerline{harris\@math.jussieu.fr \hfill ssk\@math.umd.edu}

\bye